%&amstex
\input amstex
\input amsppt.sty
\magnification=\magstep1
\hsize=32truecc
\vsize=22.2truecm
\baselineskip=16truept
\NoBlackBoxes
\TagsOnRight \pageno=1 \nologo

\def\l{\left}
\def\r{\right}
\def\bg{\bigg}
\def\({\bg(}
\def\[{\bg\lfloor}
\def\){\bg)}
\def\]{\bg\rfloor}
\def\t{\text}
\def\f{\frac}

\def\bi{\binom}
\def\eq{\equiv}

\def\ls{\leqslant}
\def\gs{\geqslant}
\def\mo{\roman{mod}}
\def\Li{\roman{Li}}

\def\da{\delta}

\def\Proof{\noindent{\it Proof}}

\def\Remark{\medskip\noindent{\it  Remark}}

\hbox {Internat. J. Math. 26(2015), no.\,8, 1550055 (23 pages).}\bigskip \topmatter
\title A new series for $\pi^3$ and related congruences\endtitle
\author Zhi-Wei Sun\endauthor
\leftheadtext{Zhi-Wei Sun} \rightheadtext{A new series for $\pi^3$
and related congruences}
\affil Department of Mathematics, Nanjing University\\
 Nanjing 210093, People's Republic of China
  \\  zwsun\@nju.edu.cn
  \\ {\tt http://math.nju.edu.cn/$\sim$zwsun}
\endaffil
\abstract Let $H_n^{(2)}$ denote the second-order harmonic number
$\sum_{0<k\ls n}1/k^2$ for $n=0,1,2,\ldots$. In this paper we obtain
the following identity:
$$\sum_{k=1}^\infty\f{2^kH_{k-1}^{(2)}}{k\bi{2k}k}=\f{\pi^3}{48}.$$
We explain how we found the series and develop related congruences
involving Bernoulli or Euler numbers; for example, it is shown that
$$\sum_{k=1}^{p-1}\f{\bi{2k}k}{2^k}H_k^{(2)}\eq-E_{p-3}\ (\mo\ p)$$
for any prime $p>3$, where $E_0,E_1,E_2,\ldots$ are Euler numbers.
Motivated by the Amdeberhan-Zeilberger identity $\sum_{k=1}^\infty(21k-8)/(k^3\bi{2k}k^3)=\pi^2/6$, we also establish the congruence
$$\sum_{k=1}^{(p-1)/2}\f{21k-8}{k^3\bi{2k}k^3}\eq(-1)^{(p+1)/2}4E_{p-3}\pmod p$$
for each prime $p>3$.
\endabstract
\thanks 2010 {\it Mathematics Subject Classification}. Primary 11B65, 11Y60;
Secondary 05A10, 05A19, 11A07, 11B68, 33E99.
\newline\indent {\it Keywords}. Series for $\pi^3$, central binomial coefficients, congruences modulo prime powers,
Bernoulli and Euler numbers.
\newline\indent This research was supported by the National Natural Science
Foundation (Grant No. 11171140) of China, and the initial version of this paper was posted to {\tt arXiv} in 2010.
\endthanks
\endtopmatter
\document

\heading{1. Introduction}\endheading

Series with summations related to $\pi$ have a long history. Leibniz
and Euler got the famous identities
$$\sum_{k=0}^\infty\f{(-1)^k}{2k+1}=\f{\pi}4\ \ \t{and}\ \
\sum_{k=1}^\infty\f1{k^2}=\f{\pi^2}6$$ respectively. Though there
exist many series for $\pi$ and $\pi^2$ (see, e.g., [Ma]), there are very few
interesting series for $\pi^3$. The most well-known series for
$\pi^3$ is the following one:
$$\sum_{k=0}^\infty\f{(-1)^k}{(2k+1)^3}=\f{\pi^3}{32}.\tag1.1$$
In 1985 Zucker [Z, (2.23)] showed that
$$\sum_{k=0}^\infty\f{\bi{2k}k}{(2k+1)^3 16^k}=\f{7\pi^3}{216}.\tag1.2$$

Recall that harmonic numbers are those rational numbers
$$H_n:=\sum_{0<k\ls n}\f1k\qquad\ (n=0,1,2,\ldots),$$
and {\it harmonic numbers of the second order} are defined by
$$H^{(2)}_n:=\sum_{0<k\ls n}\f1{k^2}\qquad (n=0,1,2,\ldots).$$

Now we give our first result which appears to be new and curious.

\proclaim{Theorem 1.1} We have the following identity:
$$\sum_{k=1}^\infty\f{2^kH_{k-1}^{(2)}}{k\bi{2k}k}=\f{\pi^3}{48}.\tag1.3$$
\endproclaim
\Remark\ 1.1. The author noted that  {\tt Mathematica 7} could not
evaluate the series in (1.3).

\medskip
By Stirling's formula
$$n!\sim\sqrt{2\pi n}\l(\f n e\r)^n\ \ (n\to+\infty)$$
and thus
$$\bi{2k}k\sim \f{4^k}{\sqrt{k\pi}}\ \ \ (k\to+\infty).$$
Note also that $H_n^{(2)}\to\zeta(2)=\pi^2/6$ as $n\to+\infty$.
Therefore
$$\f{2^kH_{k-1}^{(2)}}{k\bi{2k}k}\sim\f{\zeta(2)\sqrt{\pi}}{2^k\sqrt
k}\quad\ \ (k\to+\infty).$$ So the series in (1.3) converges much
faster than the series in (1.1) (but slower than the series in (1.2)). Using {\tt Mathematica 7}
we found that for $n\gs500$ we have
  $$\l|\f{s_n}{\pi^3/48}-1\r|<\f1{10^{150}}$$
where $s_n:= \sum_{k=1}^n2^kH_{k-1}^{(2)}/(k\bi{2k}k)$.

The reader may wonder how the author discovered (1.3) which gives a series for $\pi^3$ of a new type.
Now we present some explanations.

Let $p$ be an odd prime. In [Su3] and [Su4] the author proved the congruences
$$\sum_{k=0}^{p-1}\f{\bi{2k}k}{2^k}\eq(-1)^{(p-1)/2}-p^2E_{p-3}\ (\mo\ p^3)\tag1.4$$
and
$$\sum_{k=0}^{p-1}\bi{p-1}k\f{\bi{2k}k}{(-2)^k}\eq(-1)^{(p-1)/2}2^{p-1}\ (\mo\ p^3)\tag1.5$$
respectively, where $E_0,E_1,E_2,\ldots$ are Euler numbers given by
$E_0=1$ and the recursion
$$ \sum^n\Sb k=0\\2\mid k\endSb \bi nk E_{n-k}=0\quad\ (n=1,2,3,\ldots).$$
For $k=0,\ldots,p-1$, clearly we have
$$\align&\bi{p-1}k(-1)^k=\prod_{0<j\ls k}\l(1-\f pj\r)
\\\eq& 1-pH_k+\f{p^2}2\sum_{0<i<j\ls k}\f2{ij}
=1-pH_k+\f{p^2}2(H_k^2-H_k^{(2)})\ (\mo\ p^3).
\endalign$$
So, in view of (1.4) and (1.5), it is natural to investigate
$$\sum_{k=0}^{p-1}\f{\bi{2k}k}{2^k}H_k^{(2)}\ \mo\ p,\ \ \sum_{k=0}^{p-1}\f{\bi{2k}k}{2^k}H_k^2\ \mo\ p,
\ \ \sum_{k=0}^{p-1}\f{\bi{2k}k}{2^k}H_k\ \mo\ p^2.$$
This led the author to obtain the following result.

\proclaim{Theorem 1.2} Let $p>3$ be a prime. Then
$$\sum_{k=0}^{p-1}\f{\bi{2k}k}{2^k}H_k^{(2)}\eq-E_{p-3}\ (\mo\ p).\tag1.6$$
\endproclaim
\Remark\ 1.2.  Let $p$ be an odd prime. We are also able to show that
$$\align\sum_{k=0}^{p-1}\f{\bi{2k}k}{2^k}H_k^2\eq&\l(\f{-1}p\r)\f{q_p(2)^2}2-\f{E_{p-3}}2\ (\mo\ p),\tag1.7
\\\sum_{k=0}^{p-1}\f{\bi{2k}k}{2^k}H_k\eq&\l(\f{-1}p\r)\f{H_{(p-1)/2}}2-pE_{p-3}\ (\mo\ p^2),\tag1.8
\endalign$$
where $(\f{\cdot}p)$ denotes the Legendre symbol, and $q_p(2)$ stands for the
Fermat quotient $(2^{p-1}-1)/p$. Recall that
in 1938 Lehmer [L] proved the congruence
$$H_{(p-1)/2}\eq-2q_p(2)+p\,q_p(2)^2\ (\mo\ p^2).\tag1.9$$

\medskip

In view of certain correspondence between series for the zeta function or powers of $\pi$
and congruences involving Bernoulli or Euler numbers revealed in the authors' papers [Su2] and [Su3],
the congruence (1.6) suggests that we should consider the series
$\sum_{k=0}^\infty\bi{2k}kH_k^{(2)}/2^k$. Since this series diverges, we should seek for certain transformation.
Let $p$ be an odd prime. By [Su2, Lemma 2.1] and [T],
$$\f1p\bi{2(p-k)}{p-k}\eq-\f2{k\bi{2k}k}\ (\mo\ p)\quad\t{for}\ k=1,\ldots,\f{p-1}2.$$
Thus, if $p>3$ then
$$\align\sum_{k=0}^{p-1}\f{\bi{2k}k}{2^k}H_k^{(2)}\eq&\sum_{k=1}^{(p-1)/2}\f{k\bi{2k}k}{k2^k}H_k^{(2)}
\eq\sum_{k=1}^{(p-1)/2}\(\f{H_k^{(2)}}{k2^k}\cdot\f{-2p}{\bi{2(p-k)}{p-k}}\)
\\\eq&\sum_{p/2<k<p}\(\f{H_{p-k}^{(2)}}{(p-k)2^{p-k}}\cdot\f{-2p}{\bi{2k}k}\)
\\\eq& -p\sum_{p/2<k<p}\f{2^kH_{k-1}^{(2)}}{k\bi{2k}k}\eq-p\sum_{k=1}^{p-1}\f{2^kH_{k-1}^{(2)}}{k\bi{2k}k}\ (\mo\ p)
\endalign$$
since $2^p\eq2\ (\mo\ p)$ and
$$-H_{p-k}^{(2)}\eq H_{p-1}^{(2)}-H_{p-k}^{(2)}\eq H_{k-1}^{(2)}\ (\mo\ p).$$
Therefore the congruence in (1.6) is equivalent to
$$p\sum_{k=1}^{p-1}\f{2^kH_{k-1}^{(2)}}{k\bi{2k}k}\eq E_{p-3}\ (\mo\ p).\tag 1.6$'$ $$
Motivated by (1.6$'$) the author found (1.3).

\medskip
Now we state our third theorem which is close to Theorem 1.2.

\proclaim{Theorem 1.3} Let $p$ be an odd prime. If $p>3$, then
$$\sum_{k=1}^{p-1}\f{\bi{2k}k}{4^k}H_k\eq 2-2p+4p^2q_p(2)\ (\mo\ p^3).\tag1.10$$
We also have
$$\sum_{k=0}^{(p-1)/2}\f{\bi{2k}k}{4^k}H_k^{(2)}\eq-4q_p(2)\ (\mo\ p)\tag1.11$$
and
$$\sum_{k=1}^{(p-1)/2}\f{\bi{2k}k}{k4^k}H_k^{(2)}\eq\f{B_{p-3}}2\ (\mo\ p),\tag1.12$$
where $B_0,B_1,B_2,\ldots$ are Bernoulli numbers.
\endproclaim

In 1997 T. Amdeberhan and D. Zeilberger [AZ] obtained that
$$\sum_{k=1}^\infty\f{21k-8}{k^3\bi{2k}k^3}=\zeta(2)=\f{\pi^2}6.$$
We are able to establish the following result
related to the Amdeberhan-Zeilberger series.

\proclaim{Theorem 1.4} Let $p>3$ be a prime. Then
$$\sum_{k=0}^{(p-1)/2}(21k+8)\bi{2k}k^3\eq 8p+(-1)^{(p-1)/2}32p^3E_{p-3}\ (\mo\ p^4)\tag1.13$$
and hence
$$\sum_{k=1}^{(p-1)/2}\f{21k-8}{k^3\bi{2k}k^3}\eq(-1)^{(p+1)/2}4E_{p-3}\pmod p.\tag1.14$$
\endproclaim
\Remark\ 1.3. In [Su3] the author showed that
$$\sum_{k=0}^{p-1}(21k+8)\bi{2k}k^3\eq 8p+16p^4B_{p-3}\ (\mo\ p^5)\tag1.15$$
for any odd prime $p$. However, (1.13) is much more sophisticated than this congruence involving $B_{p-3}$.

\medskip

The next section is devoted to the proof of Theorem 1.1.
We are going to show Theorems 1.2--1.3 and Theorem 1.4 in Sections 3 and 4 respectively.
Section 5 contains some conjectures of the author for further research.

\heading{2. Proof of Theorem 1.1}\endheading

Set $$S:=\sum_{k=1}^\infty\f{2^kH_{k-1}^{(2)}}{k\bi{2k}k}.$$
Then
$$S=\sum_{k=0}^\infty\f{2^{k+1}H_k^{(2)}}{(k+1)\bi{2k+2}{k+1}}
=\sum_{k=0}^\infty\f{2^{k}H_k^{(2)}}{(k+1)\bi{2k+1}k}
=\sum_{k=0}^\infty\f{2^{k}H_k^{(2)}\Gamma(k+1)^2}{\Gamma(2k+2)}.$$
Recall the well-known fact that
$$B(a,b):=\int_0^1x^{a-1}(1-x)^{b-1}dx=\f{\Gamma(a)\Gamma(b)}{\Gamma(a+b)}\
\quad \t{for any}\ a,b>0.$$
So we have
$$\align S=&\sum_{k=0}^\infty2^{k}H_k^{(2)}\int_0^1x^k(1-x)^kdx
=\sum_{k=0}^\infty\f{H_k^{(2)}}{2^k}\int_0^1(1-(2x-1)^2)^kdx
\\=&\sum_{k=0}^\infty\f{H_k^{(2)}}{2^{k+1}}\int_{-1}^1(1-t^2)^kdt
=\sum_{k=1}^\infty \f{H_k^{(2)}}{2^{k}}\int_0^1(1-t^2)^kdt.
\endalign$$

 Observe that if $0\ls t\ls 1$ then
$$\align \sum_{k=1}^\infty H_k^{(2)}\l(\f{1-t^2}2\r)^k=&\sum_{k=1}^\infty\sum_{j=1}^k\f1{j^2}\l(\f{1-t^2}2\r)^k
=\sum_{j=1}^\infty\f1{j^2}\sum_{k=j}^\infty\l(\f{1-t^2}2\r)^k
\\=&\sum_{j=1}^\infty\f1{j^2}\l(\f{1-t^2}2\r)^j\f1{1-(1-t^2)/2}
\\=&\f2{1+t^2}\Li_2\l(\f{1-t^2}2\r),
\endalign$$
where the dilogarithm $\Li_2(x)$ is given by
$$\Li_2(x):=\sum_{n=1}^\infty\f{x^n}{n^2}\qquad \ (|x|<1).$$
Therefore
$$\f S2=\int_0^1\f1{1+t^2}\Li_2\l(\f{1-t^2}2\r)dt=\int_0^1\Li_2\l(\f{1-t^2}2\r)(\arctan t)'dt.$$

 Note that
$$\Li_2'(x)=\sum_{n=1}^\infty\f{x^{n-1}}{n}=-\f{\log(1-x)}x$$
and hence
$$\f{d}{dt}\Li_2\l(\f{1-t^2}2\r)=-\f{\log(1-(1-t^2)/2)}{(1-t^2)/2}\times(-t)=\f{2t}{1-t^2}\log\f{1+t^2}2.$$
Thus, using integration by parts we obtain
$$\align \f S2=&\Li_2\l(\f{1-t^2}2\r)\arctan t\ \bigg|_{t=0}^1-\int_0^1(\arctan t)\f{2t}{1-t^2}\log\f{1+t^2}2dt
\\=&\int_0^1(\arctan t)\l(\f1{1+t}-\f1{1-t}\r)\log\f{1+t^2}2dt
\\=&\int_0^1\f{\arctan t}{1+t}\log\f{1+t^2}2dt-\int_0^{-1}\f{\arctan t}{1+t}\log\f{1+t^2}2dt
\\=&\int_{-1}^1\f{\arctan t}{1+t}\log\f{1+t^2}2dt.
\endalign$$

 Finally, inputting the Mathematica command
\newline
\newline
\indent\indent\indent {\tt Integrate[ArcTan[t]Log[(1+t${}^\wedge$2)/2]/(1+t),\{t,-1,1\}]}
\newline
\newline
we then obtain from {\tt Mathematica 7} that
$$\int_{-1}^1\f{\arctan t}{1+t}\log\f{1+t^2}2dt=\f{\pi^3}{96}.$$
Thus $S=\pi^3/48$ as desired. We are done.

\heading{3. Proofs of Theorems 1.2 and 1.3}\endheading

 We first state some basic facts which will be used very often.
 For any prime $p>3$ we have
 $$\sum_{k=1}^{(p-1)/2}\f1{k^2}\eq\f12\sum_{k=1}^{p-1}\f1{k^2}\eq0\ (\mo\ p)$$
 since $\sum_{j=1}^{p-1}(2j)^{-2}\eq\sum_{k=1}^{p-1}k^{-2}\ (\mo\ p)$.  If $p$ is an odd prime, then
 $$\bi {(p-1)/2}k\eq\bi{-1/2}k=\f{\bi{2k}k}{(-4)^k}\pmod p\quad\t{for all}\ k=0,\ldots,p-1.\tag3.1$$
For any $n=0,1,2,\ldots$ we have the identity
 $$\sum_{k=0}^n(-1)^k\bi xk=(-1)^n\bi{x-1}n\tag3.2$$
 which can be found in [G, (1.5)].

\proclaim{Lemma 3.1} For any positive integer $n$, we have the identities
$$\sum_{k=1}^n\bi nk\f{(-1)^{k-1}}k=H_n\tag3.3$$
and
$$\sum_{k=1}^n\bi nk\f{(-1)^{k-1}}kH_k=H_n^{(2)}.\tag3.4$$
\endproclaim
\Proof. (3.3) and (3.4) follow from [G, (1.45)] and an identity of V. Hern\'andez [He]
respectively. Below we give a simple proof of (3.4). In view of the binomial inversion formula (cf. (5.48) of [GKP, pp.\,192-193]),
(3.4) holds for all $n=1,2,3,\ldots$
if and only if for any positive integer $n$ we have
$$\sum_{k=1}^n\bi nk(-1)^k H_k^{(2)}=-\f{H_n}n.\tag3.4$'$ $$
In fact, in view of (3.2) and (3.3), we get
$$\align\sum_{k=1}^n\bi{n}k(-1)^k\sum_{j=1}^k\f1{j^2}=&\sum_{j=1}^n\f1{j^2}\(\sum_{k=0}^n\bi nk(-1)^k-\sum_{k=0}^{j-1}\bi nk(-1)^k\)
\\=&\sum_{j=1}^n\f{(-1)^j}{j^2}\bi{n-1}{j-1}=\f1n\sum_{j=1}^n\f{(-1)^j}j\bi nj=-\f{H_n}n
\endalign$$
and hence (3.4$'$) holds. \qed

\proclaim{Lemma 3.2} Let $p=2n+1$ be an odd prime and let $m$ be an integer with $m\not\eq0,4\ (\mo\ p)$.
Then
$$\sum_{k=1}^n\f{\bi{2k}k}{m^k}H_k^{(2)}\eq-\l(\f{m(m-4)}p\r)\sum_{k=1}^n\f{\bi{2k}kH_k}{k(4-m)^k}\ (\mo\ p).\tag3.5$$
In particular,
$$\sum_{k=1}^{(p-1)/2}\f{\bi{2k}k}{2^k}H_k^{(2)}\eq-\l(\f{-1}p\r)\sum_{k=1}^{(p-1)/2}\f{\bi{2k}kH_k}{k2^k}\ (\mo\ p).\tag3.6$$
\endproclaim
\Proof. Clearly it suffices to prove (3.5).

In view of (3.4), we have
$$\align \sum_{k=1}^n\bi nk\l(-\f 4m\r)^k H_k^{(2)}=&\sum_{k=1}^n\f{\bi{n}k(-4)^k}{m^k}\sum_{j=1}^k\bi kj\f{(-1)^{j-1}}jH_j
\\=&\sum_{j=1}^n\f{(-1)^{j-1}}jH_j\sum_{k=j}^n\bi nk\bi kj\l(-\f 4m\r)^k
\\=&\sum_{j=1}^n\f{(-1)^{j-1}}jH_j\bi nj\sum_{k=j}^n\bi {n-j}{k-j}\l(-\f 4m\r)^k
\\=&\sum_{j=1}^n\bi nj\f{(-1)^{j-1}}jH_j\l(-\f 4m\r)^j\l(1-\f 4m\r)^{n-j}
\\=&-\f1{m^n}\sum_{j=1}^n\bi nj\f{4^jH_j}j(m-4)^{n-j}.
\endalign$$
So, with the help of (3.1), we obtain
$$\sum_{k=1}^n\f{\bi{2k}k}{m^k}H_k^{(2)}\eq-\l(\f{m(m-4)}p\r)\sum_{j=1}^n\f{\bi{2j}j(-1)^jH_j}{j(m-4)^j}\pmod{p}.$$
This proves (3.5). We are done. \qed

\proclaim{Lemma 3.3} Let $n$ be any positive integer. Then
$$\sum_{k=1}^n\bi nk\f{(-2)^k}kH_k=-2\sum^n\Sb k=1\\2\nmid k\endSb\f{H_n-H_{n-k}}k.\tag3.7$$
\endproclaim
\Proof. Let $S_n$ denote the left-hand side of (3.7). Observe that
$$\align S_n=&\sum_{k=1}^n\bi nk\f{(-2)^k}k\sum_{j=1}^k\int_0^1x^{j-1}dx=\int_0^1\sum_{k=1}^n\bi nk\f{(-2)^k}k\cdot\f{x^k-1}{x-1} dx
\\=&\int_0^1\int_0^1\sum_{k=0}^n\bi nk\f{(-2x)^k-(-2)^k}{x-1}y^{k-1}dydx
\\=&\int_0^1\int_0^1\f{(1-2xy)^n-(1-2y)^n}{(x-1)y}dydx
\\=&-2\int_0^1\int_0^1\sum_{k=1}^{n}(1-2xy)^{k-1}(1-2y)^{n-k}dxdy.
\endalign$$
Clearly,
$$\align\int_0^1(1-2xy)^{k-1}dx=\f{(1-2xy)^{k}}{-2ky}\bigg|_{x=0}^1
=\f{(1-2y)^{k}-1}{k(1-2y-1)}=\f1{k}\sum_{j=1}^{k}(1-2y)^{j-1}.
\endalign$$
Therefore
$$\align S_n=&-2\int_0^1\sum_{k=1}^n\f1k\sum_{j=1}^k(1-2y)^{n-k+j-1}dy=\sum_{1\ls j\ls k\ls n}\f{(1-2y)^{n-k+j}}{k(n-k+j)}\bigg|_{y=0}^1
\\=&\sum_{1\ls j\ls k\ls n}\f{(-1)^{n-k+j}-1}{k(n-k+j)}=\sum_{i=1}^n\f{(-1)^i-1}i\sum_{j=1}^i\f1{n+j-i}
\\=&-2\sum^n\Sb i=1\\2\nmid i\endSb\f1i(H_n-H_{n-i}).
\endalign$$
This completes the proof of (3.7). \qed

\proclaim{Lemma 3.4} Let $p>3$ be a prime. Then
$$\bi{p-1}{(p-1)/2}\eq(-1)^{(p-1)/2}4^{p-1}\ (\mo\ p^3).\tag3.8$$
\endproclaim
\Remark\ 3.1. (3.8) is a famous congruence of Morley [Mo].
\medskip

\proclaim{Lemma 3.5} Let $p>3$ be a prime. Then
$$\sum^{(p-1)/2}\Sb k=1\\2\nmid k\endSb\f{H_k}k\eq\f34q_p(2)^2-\l(\f{-1}p\r)\f{E_{p-3}}2\ (\mo\ p)\tag3.9$$
and
$$\sum^{(p-1)/2}\Sb k=1\\2\mid k\endSb\f{H_k}k\eq\f54q_p(2)^2+\l(\f{-1}p\r)\f{E_{p-3}}2\ (\mo\ p).\tag3.10$$
\endproclaim
\Proof.  Set $n=(p-1)/2$. Clearly it suffices to show that
$$\sum_{k=1}^n\f{H_k}k\eq2q_p(2)^2\ (\mo\ p)\tag3.11$$ and
$$\sum_{k=1}^n\f{(-1)^k}kH_k\eq\f{q_p(2)^2}2+\l(\f{-1}p\r)E_{p-3}\ (\mo\ p).\tag3.12$$

Let $\da\in\{0,1\}$. For $r=0,\ldots,p-1$ we obviously have
$$(-1)^r\bi{p-1}{r}=\prod_{0<s\ls r}\l(1-\f ps\r)\eq1-pH_r\ (\mo\ p^2).$$
Thus
$$\align p\sum_{k=1}^n\f{(-1)^{\da k}}kH_{k-1}\eq&\sum_{k=1}^n\f{(-1)^{\da k}}k\l(1-(-1)^{k-1}\bi{p-1}{k-1}\r)
\\=&\sum_{k=1}^n\f{(-1)^{\da k}}k+\f1p\sum_{k=1}^n(-1)^{(\da+1)k}\bi pk\ (\mo\ p^2)
\endalign$$
and hence
$$p\sum_{k=1}^n\f{(-1)^{\da k}}kH_{k}\eq
\sum_{k=1}^n(-1)^{\da k}\l(\f1k+\f p{k^2}\r)+\f1p\sum_{k=1}^n(-1)^{(\da+1)k}\bi pk\ (\mo\ p^2).\tag3.13$$

Putting $\da=0$ in (3.13) and recalling (3.2) and the congruence $\sum_{k=1}^n1/k^2\eq0\ (\mo\ p)$, we get
$$p\sum_{k=1}^n\f{H_k}k\eq H_n+\f{(-1)^n\bi{p-1}{n}-1}p\pmod{p^2}.$$
With the helps of (1.9) and (3.8), we have
$$p\sum_{k=1}^n\f{H_k}k\eq-2q_p(2)+pq_p(2)^2+\f{(1+p\,q_p(2))^2-1}p\eq 2p\,q_p(2)^2 \ (\mo\ p^2)$$
which yields (3.11).
Taking $\da=1$ in (3.13) and using the congruence $\sum_{k=1}^n1/k^2\eq0\ (\mo\ p)$, we obtain
$$\align p\sum_{k=1}^n\f{(-1)^k}kH_k\eq&\sum_{k=1}^n\f{(-1)^k+1}k+p\sum_{k=1}^n\f{(-1)^k+1}{k^2}-H_n+\f{2^p-2}{2p}
\\=&H_{\lfloor p/4\rfloor}+\f p2\sum_{j=1}^{\lfloor p/4\rfloor}\f1{j^2}-H_n+q_p(2)\pmod{p^2}.
\endalign$$
Let's recall (1.9) and note that
$$\sum_{k=1}^{\lfloor p/4\rfloor}\f1{k^2}\eq4\l(\f{-1}p\r)E_{p-3}\ (\mo\ p)\tag3.14$$
and
$$H_{\lfloor p/4\rfloor}\eq-3q_p(2)+\f32p\,q_p(2)^2-\l(\f{-1}p\r)pE_{p-3}\ (\mo\ p^2)$$
by Lehmer [L, (20)] and [S2, Corollary 3.3] respectively. Therefore,
$$\align p\sum_{k=1}^n\f{(-1)^k}kH_k\eq&-3q_p(2)+\f32p\,q_p(2)^2-\l(\f{-1}p\r)pE_{p-3}
\\&+2p\l(\f{-1}p\r)E_{p-3}+(2q_p(2)-p\,q_p(2)^2)+q_p(2)
\\=&\f p2 q_p(2)^2+\l(\f{-1}p\r)pE_{p-3}\pmod{p^2}
\endalign$$
and hence (3.12) holds. We are done. \qed

\proclaim{Lemma 3.6} Let $p$ be an odd prime. Then
$$\sum^{p-1}\Sb k=1\\4\mid k-2\endSb\f{H_k}k\eq\f 3{16}q_p(2)^2\pmod{p}.\tag3.15$$
If $p>3$, then we also have
$$\sum^{p-1}\Sb k=1\\4\mid k\endSb\f{H_k}k\eq\f 5{16}q_p(2)^2\pmod{p}.\tag3.16$$
\endproclaim
\Proof. As $H_{p-k}=H_{p-1}-\sum_{0<j<k}1/(p-j)\eq H_{k-1}\ (\mo\ p)$ for $k=1,\ldots,p-1$, we have
$$\align p\sum^{p-1}\Sb k=1\\4\mid k-2\endSb\f{H_k}k=&p\sum^{p-1}\Sb k=1\\4\mid k-p+2\endSb\f{H_{p-k}}{p-k}
\\\eq&-\sum^{p-1}\Sb k=1\\4\mid k+p\endSb\f{pH_{k-1}}k
\eq\sum^{p-1}\Sb k=1\\4\mid k+p\endSb\f{(-1)^{k-1}\bi{p-1}{k-1}-1}k
\\=&\sum^{p-1}\Sb k=1\\4\mid k+p\endSb\f1k\bi{p-1}{k-1}-\sum^{p-1}\Sb k=1\\4\mid k+p\endSb\f1k\pmod{p^2}.
\endalign$$
Note that
$$2\sum^{p-1}\Sb k=1\\4\mid k+p\endSb\f1k\bi{p-1}{k-1}=q_p(2)-\f{(\f 2p)2^{(p-1)/2}-1}p=\f{2^{p-1}-(\f2p)2^{(p-1)/2}}p$$
by [Su1, Corollary 3.1] and that
$$\sum^{p-1}\Sb k=1\\4\mid k+p\endSb\f1k\eq\f{q_p(2)}4-\f p8 q_p(2)^2\pmod{p^2}$$
by [S2, Corollary 3.1]. Therefore
$$\align p\sum^{p-1}\Sb k=1\\4\mid k-2\endSb\f{H_k}k\eq&\f{2^{p-1}-(\f2p)2^{(p-1)/2}}{2p}-\f{2^{p-1}-1}{4p}+\f p8q_p(2)^2
\\=&p\l(\f{(\f2p)2^{(p-1)/2}-1}{2p}\r)^2+\f p8q_p(2)^2
\\\eq& p\l(\f{2^{p-1}-1}{4p}\r)^2+\f p8q_p(2)^2=\f3{16}p\,q_p(2)^2
\pmod{p^2}\endalign$$
and hence (3.15) follows. When $p>3$ we can prove (3.16) in a similar way. \qed

\medskip
\noindent{\it Proof of Theorem 1.2}. Set $n=(p-1)/2$. In view of (3.1) and (3.6), it suffices to show
$$\sum_{k=1}^n\bi nk\f{(-2)^k}kH_k\eq\l(\f{-1}p\r)E_{p-3}\pmod{p}.$$
For each $k=1,\ldots,n$, evidently
$$H_n-H_{n-k}=\sum_{j=0}^{k-1}\f1{n-j}\eq-2\sum_{j=0}^{k-1}\f1{2j+1}=-2\l(H_{2k}-\f{H_k}2\r)\ (\mo\ p).$$
Thus, in light of (3.7), (3.15) and (3.9), we have
$$\align \sum_{k=1}^n\bi nk\f{(-2)^k}kH_k\eq& 4\sum^n\Sb k=1\\2\nmid k\endSb\l(\f{H_{2k}}k-\f{H_k}{2k}\r)
=8\sum^{p-1}\Sb j=1\\4\mid j-2\endSb\f{H_j}j-2\sum^n\Sb k=1\\2\nmid k\endSb\f{H_k}k
\\\eq&\f 32q_p(2)^2-\f32q_p(2)^2+\l(\f{-1}p\r)E_{p-3}\pmod{p}
\endalign$$
as desired. This concludes the proof. \qed

\proclaim{Lemma 3.7} Let $p>3$ be a prime. Then
$$\sum_{k=1}^{p-1}\f{H_k}{k^2}\eq B_{p-3}\pmod{p}\tag3.17$$
and
$$\sum_{k=1}^{(p-1)/2}\f1{k^3}\eq-2B_{p-3}\pmod{p}.\tag3.18$$
\endproclaim
\Remark\ 3.2. (3.17) appeared as [ST, (5.4)], and (3.18) follows from [S1, Corollary 5.2(b)].

\proclaim{Lemma 3.8} For any positive integer $m$ and nonnegative integer $n$ we have
$$\sum_{k=0}^n\bi nk\f{(-1)^k}{k+m}=\f1{m\bi{m+n}m}.\tag3.19$$
\endproclaim
\Remark\ 3.3. (3.19) can be found in [G, (1.43)].

\medskip\noindent{\it Proof of Theorem 1.3}. Observe that
$$\align\sum_{k=1}^{p-1}\f{\bi{2k}k}{4^k}H_k=&\sum_{k=1}^{p-1}\bi{-1/2}k(-1)^k\sum_{j=1}^k\f1j
\\=&\sum_{j=1}^{p-1}\f1j\(\sum_{k=0}^{p-1}\bi{-1/2}k(-1)^k-\sum_{k=0}^{j-1}\bi{-1/2}k(-1)^k\).
\endalign$$
Applying (3.2) we get
$$\align\sum_{k=1}^{p-1}\f{\bi{2k}k}{4^k}H_k=&\sum_{j=1}^{p-1}\f1j\((-1)^{p-1}\bi{-1/2-1}{p-1}-(-1)^{j-1}\bi{-1/2-1}{j-1}\)
\\=&\bi{-3/2}{p-1}H_{p-1}-2\sum_{j=1}^{p-1}(-1)^j\f{-1/2}j\bi{-3/2}{j-1}
\\=&\bi{-3/2}{p-1}H_{p-1}-2\(\sum_{j=0}^{p-1}(-1)^j\bi{-1/2}j-1\)
\\=&\bi{-3/2}{p-1}H_{p-1}-2\bi{-1/2-1}{p-1}+2.
\endalign$$
Now assume $p>3$. Note that
$$\bi{-3/2}{p-1}=\f p{-1/2}\bi{-1/2}p=-2p\f{\bi{2p}p}{(-4)^p}=p\f{\bi{2p-1}{p-1}}{4^{p-1}}
\eq\f p{4^{p-1}}\ (\mo\ p^4)$$
since $\bi{2p-1}{p-1}\eq1\ (\mo\ p^3)$ by Wolstenholme's theorem (see, e.g., [HT]).
In view of Wolstenholme's congruence $H_{p-1}\eq0\ (\mo\ p^2)$,
by the above we have
$$\align\sum_{k=1}^{p-1}\f{\bi{2k}k}{4^k}H_k\eq&2-2\bi{-3/2}{p-1}\eq2-\f {2p}{(1+p\,q_p(2))^2}
\\\eq&2- 2p(1-p\,q_p(2))^2=2-2p+4p^2q_p(2)\ (\mo\ p^3).
\endalign$$
So (1.10) holds.

Below we write $p=2n+1$.
Combining (3.1), (3.4$'$) and (1.9), we get
$$\sum_{k=0}^n\f{\bi{2k}k}{4^k}H_k^{(2)}\eq-\f{H_{(p-1)/2}}{(p-1)/2}\eq2H_{(p-1)/2}\eq-4q_p(2)\ (\mo\ p).$$
This proves (1.11).

In view of (3.4) and (3.19), we have
$$\align\sum_{k=1}^n\bi nk\f{(-1)^k}kH_k^{(2)}=&\sum_{k=1}^n\bi nk\f{(-1)^k}k\sum_{j=1}^k\bi kj(-1)^{j-1}\f{H_j}j
\\=&-\sum_{j=1}^n\f{H_j}j\bi nj\sum_{k=j}^n\f{(-1)^{k-j}}k\bi{n-j}{k-j}
\\=&-\sum_{j=1}^n\f{H_j}j\bi nj\f1{j\bi nj}=-\sum_{j=1}^n\f{H_j}{j^2}.
\endalign$$
Observe that
$$\align\sum_{k=1}^{p-1}\f{H_k}{k^2}=&\sum_{k=1}^n\(\f{H_k}{k^2}+\f{H_{p-k}}{(p-k)^2}\)
\\\eq&\sum_{k=1}^n\(\f{H_k}{k^2}+\f{H_{k-1}}{k^2}\)=2\sum_{k=1}^n\f{H_k}{k^2}-\sum_{k=1}^n\f1{k^3}\pmod{p}.
\endalign$$
Therefore, with the help of (3.1) we have
$$\sum_{k=1}^n\f{\bi{2k}k}{4^k}H_k^{(2)}\eq-\sum_{k=1}^n\f{H_k}{k^2}
\eq-\f12\(\sum_{k=1}^{p-1}\f{H_k}{k^2}+\sum_{k=1}^n\f1{k^3}\)\pmod{p}.$$
Now applying Lemma 3.7 we immediately get the desired (1.12).

The proof of Theorem 1.3 is now complete. \qed

\heading{4. Proof of Theorem 1.4}\endheading

\proclaim{Lemma 4.1} For any positive integer $n$, we have the following identities:
$$\align
\\\sum_{k=0}^n\bi nk\bi{n+k}k(-1)^kH_k=&2(-1)^nH_n,\tag4.1
\\\sum_{k=0}^n\bi nk\bi{n+k}k(-1)^kH_k^{(2)}=&2(-1)^{n-1}\sum_{k=1}^n\f{(-1)^k}{k^2}.\tag4.2
\endalign$$
\endproclaim
\Remark\ 4.1. (4.1) and (4.2) can be found in [OS] and [Pr].

\proclaim{Lemma 4.2} Let $p=2n+1$ be an odd prime, and let $k\in\{0,\ldots,n\}$. Then
$$\aligned\f{\bi {n+k}k}{\bi{2k}k/4^k}\eq&1+p\sum_{j=1}^k\f1{2j-1}+\f{p^2}2\(\sum_{j=1}^k\f1{2j-1}\)^2
\\&-\f{p^2}2\sum_{j=1}^k\f1{(2j-1)^2}
\pmod{p^3}\endaligned\tag4.3$$
and
$$\aligned\f{\bi nk}{\bi{2k}k/(-4)^k}\eq&1-p\sum_{j=1}^k\f1{2j-1}+\f{p^2}2\(\sum_{j=1}^k\f1{2j-1}\)^2
\\&-\f{p^2}2\sum_{j=1}^k\f1{(2j-1)^2}
\pmod{p^3}.
\endaligned\tag4.4$$
Consequently,
$$\bi nk\bi{n+k}k(-1)^k\eq\f{\bi{2k}k^2}{16^k}\pmod{p^2}.\tag4.5$$
\endproclaim
\Proof. Observe that
$$\align\f{\bi{n+k}k}{\bi{2k}k/4^k}=&\prod_{j=1}^k\f{(n+j)/j}{(2j-1)/(2j)}=\prod_{j=1}^k\l(1+\f p{2j-1}\r)
\\\eq&1+p\sum_{j=1}^k\f1{2j-1}+\f{p^2}2 S_k\pmod{p^3},\endalign$$
where
$$S_k:=2\sum_{1\ls i<j\ls k}\f1{(2i-1)(2j-1)}=\(\sum_{j=1}^k\f1{2j-1}\)^2-\sum_{j=1}^k\f1{(2j-1)^2}.$$
This proves (4.3). Similarly,
$$\f{(-1)^k\bi nk}{\bi{2k}k/4^k}=\prod_{j=1}^k\l(1-\f p{2j-1}\r)\eq 1-p\sum_{j=1}^k\f1{2j-1}+\f{p^2}2S_k\pmod{p^3}$$
and hence (4.4) holds. Clearly (4.5) follows from (4.3) and (4.4). We are done. \qed
\Remark\ 4.2. The congruence (4.5) was first observed by van Hamme [vH].

\proclaim{Lemma 4.3} For any nonnegative integer $n$ we have
$$\sum_{k=0}^n\bi nk^2=\bi{2n}n\tag4.6$$
and
$$\sum_{k=0}^n\bi nk\f{\bi{2k}k}{(-4)^k}=\f{\bi{2n}n}{4^n}.\tag4.7$$
\endproclaim
\Remark\ 4.3. As $\bi nk=\bi n{n-k}$ and $\bi{2k}k/(-4)^k=\bi{-1/2}k$ for all $k=0,\ldots,n$,
both (4.6) and (4.7) are special cases of the well-known Chu-Vandermonde identity
$\sum_{k=0}^n\bi xk\bi y{n-k}=\bi{x+y}n$ (cf. [G, (3.1)] or (5.22) of [GKP, p.\,169]). \qed

\proclaim{Lemma 4.5}
Let $n$ be any positive integer.
Then $$t_n:=\f1{4n\bi{2n}n}\sum_{k=0}^{n-1}(21k+8)\bi{2k}k^3$$
coincides with
$$t_n':=\sum_{k=0}^{n-1}\bi{n+k-1}k^2.$$
\endproclaim
\Remark\ 4.4. In Feb. 2010, the author conjectured that $t_n$ is always an integer
and later this was confirmed by Kasper Andersen by getting $t_n=t_n'$ via the Zeilberger algorithm
(cf. [Su3, Lemma 4.1]).
\medskip

Now we are ready to prove the following auxiliary result.

\proclaim{Theorem 4.1} Let $p>3$ be a prime. Then
$$\align \sum_{k=0}^{(p-1)/2}\bi{2k}k^2\f{H_k}{16^k}\eq& 2\l(\f{-1}p\r)H_{(p-1)/2}\pmod{p^2},\tag4.8
\\\sum_{k=1}^{(p-1)/2}\bi{2k}k^2\f{H_k^{(2)}}{16^k}\eq&-4E_{p-3}\pmod{p},\tag4.9
\\\sum_{k=1}^{(p-1)/2}\bi{2k}k^2\f{H_k}{k16^k}\eq& 4\l(\f{-1}p\r)E_{p-3}\pmod{p},\tag4.10
\\\sum_{k=0}^{(p-1)/2}\bi{2k}k^2\f{H_{2k}}{16^k}\eq& \l(\f{-1}p\r)\f32H_{(p-1)/2}+pE_{p-3}\pmod{p^2}.\tag4.11
\endalign$$
\endproclaim
\Proof.  Set $n=(p-1)/2$. In view of (4.5), (4.1) implies (4.8), and (4.2) yields that
$$\sum_{k=0}^n\bi{2k}k^2\f{H_k^{(2)}}{16^k}\eq 2(-1)^{n-1}\sum_{k=1}^n\f{(-1)^k}{k^2}\pmod{p^2}.$$
Since $\sum_{k=1}^n1/k^2\eq0\ (\mo\ p)$, we have
$$\sum_{k=1}^n\f{(-1)^k}{k^2}\eq\sum_{k=1}^n\f{(-1)^k+1}{k^2}
=\f12\sum_{j=1}^{\lfloor p/4\rfloor}\f1{j^2}\eq2(-1)^nE_{p-3}\pmod{p}$$
by applying (3.14) in the last step. Now it is clear that (4.9) holds.

Next we deduce (4.10). With the helps of (3.4) and the Chu-Vandermonde identity,
we get
$$\align&\sum_{k=0}^n\bi nk\bi{n+k}k(-1)^kH_k^{(2)}
\\=&\sum_{k=1}^n\bi nk\bi{n+k}k(-1)^k\sum_{j=1}^k\bi kj\f{(-1)^{j-1}}jH_j
\\=&\sum_{j=1}^n\bi nj\f{(-1)^{j-1}}jH_j\sum_{k=j}^n\bi{n+k}k(-1)^k\bi{n-j}{k-j}
\\=&\sum_{j=1}^n\bi nj\f{(-1)^{j-1}}jH_j\sum_{k=0}^n\bi{-n-1}k\bi{n-j}{n-k}
\\=&\sum_{j=1}^n\bi nj\f{(-1)^{j-1}}jH_j\bi{-j-1}n=(-1)^{n-1}\sum_{j=1}^n\bi nj\bi{n+j}j\f{(-1)^j}jH_j.
\endalign$$
Thus, by applying (4.5) we obtain (4.10) from (4.9).

Since
$$\align \sum_{k=0}^n\bi nk^2H_{2k}^{(2)}=&\sum_{k=0}^n\bi nk^2H^{(2)}_{2(n-k)}=\sum_{k=0}^n\bi nk^2H_{p-1-2k}^{(2)}
\\\eq&-\sum_{k=0}^n\bi nk^2H_{2k}^{(2)}\ (\mo\ p),
\endalign$$
by (3.1) we have
$$\sum_{k=0}^n\f{\bi{2k}k^2}{16^k}H_{2k}^{(2)}\eq\sum_{k=0}^n\bi nk^2H_{2k}^{(2)}\eq0\ (\mo\ p)$$
and hence
$$\align\sum_{k=0}^n\f{\bi{2k}k^2}{16^k}\sum_{j=1}^k\f1{(2j-1)^2}
=&\sum_{k=0}^n\f{\bi{2k}k^2}{16^k}\l(H_{2k}^{(2)}-\f{H_k^{(2)}}4\r)
\\\eq&-\f14\sum_{k=0}^n\f{\bi{2k}k^2}{16^k}H_k^{(2)}\ (\mo\ p).
\endalign$$
Thus (4.9) implies that
$$\sum_{k=0}^n\f{\bi{2k}k^2}{16^k}\sum_{j=1}^k\f1{(2j-1)^2}\eq E_{p-3}\pmod{p}.\tag4.12$$

By [Su3, (1.7)],
$$\sum_{k=0}^n\f{\bi{2k}k^2}{16^k}\eq(-1)^n+p^2E_{p-3}\pmod{p^3},\tag4.13$$
Combining this with (4.12), we see that
$$\sum_{k=0}^n\f{\bi{2k}k^2}{16^k}\(1-p^2\sum_{j=1}^k\f1{(2j-1)^2}\)\eq(-1)^n\pmod{p^3}.$$

By (4.6) and (4.7), we have
$$\align\l(1-\f2{4^n}\r)\bi{2n}n=&\sum_{k=0}^n\bi nk\(\bi nk-\f{2\bi{2k}k}{(-4)^k}\)
\\=&\sum_{k=0}^n\f{\bi{2k}k^2}{16^k}\cdot\f{\bi nk}{\bi{2k}k/(-4)^k}\(\f{\bi nk}{\bi{2k}k/(-4)^k}-2\).
\endalign$$
Combining this with (4.4) we get
$$\l(1-\f2{4^n}\r)\bi{2n}n\eq\sum_{k=1}^n\f{\bi{2k}k^2}{16^k}\(p^2\(\sum_{j=1}^k\f1{2j-1}\)^2-1\)\pmod{p^3}.$$
By Morley's congruence (3.8),
$$\l(1-\f2{4^n}\r)\bi{2n}n+(-1)^n\eq(-1)^n(4^{2n}-2\cdot 4^n+1)=(-1)^np^2q_p(2)^2\ (\mo\ p^3).$$
Thus, in light of (4.13) we obtain
$$\sum_{k=1}^{(p-1)/2}\f{\bi{2k}k^2}{16^k}\(\sum_{j=1}^k\f1{2j-1}\)^2\eq E_{p-3}+\l(\f{-1}p\r)q_p(2)^2\pmod{p}.
\tag4.14$$

By (4.7), (4.4), (4.12) and (4.14),
$$\align&\f{\bi{2n}n}{4^n}-\sum_{k=0}^n\f{\bi{2k}k^2}{16^k}\(1-p\sum_{j=1}^k\f1{2j-1}\)
\\\eq&\f{p^2}2\sum_{k=0}^n\f{\bi{2k}k^2}{16^k}\(\(\sum_{j=1}^k\f1{2j-1}\)^2-\sum_{j=1}^k\f1{(2j-1)^2}\)
\\\eq&\f{p^2}2(-1)^nq_p(2)^2\pmod{p^3}.
\endalign$$
Combining this with (3.8) and (4.13) we obtain
$$\sum_{k=1}^n\f{\bi{2k}k^2}{16^k}\sum_{j=1}^k\f1{2j-1}\eq(-1)^n\l(-q_p(2)+\f p2q_p(2)^2\r)+pE_{p-3}\pmod{p^2}.
\tag4.15$$
Therefore, in view of (4.8) and (1.9), we have
$$\align\sum_{k=0}^n\f{\bi{2k}k^2}{16^k}H_{2k}
=&\sum_{k=1}^n\f{\bi{2k}k^2}{16^k}\(\sum_{j=1}^k\f1{2j-1}+\f{H_k}2\)
\\\eq&(-1)^n\l(-q_p(2)+\f p2q_p(2)^2\r)+pE_{p-3}+(-1)^nH_n
\\\eq&(-1)^n\f 32H_n+pE_{p-3}\ (\mo\ p^2).
\endalign$$
This proves (4.11).

So far we have completed the proof of Theorem 4.1. \qed

\medskip
\noindent{\it Proof of Theorem 1.4}. Write $p=2n+1$. Clearly
$$4(n+1)\bi{2(n+1)}{n+1}=8p\bi{2n}n\eq 8p(-1)^n4^{p-1}\pmod{p^4}$$
by Morley's congruence (3.8), and
$$\align &4^{1-p}=\l(\f1{1+p\,q_p(2)}\r)^2
\\\eq&(1-p\,q_p(2)+p^2q_p(2)^2)^2\eq 1-2p\,q_p(2)+3p^2q_p(2)^2\ (\mo\ p^3).
\endalign$$
Thus, in view of Lemma 4.5, (1.13) is reduced to
$$\aligned\sum_{k=0}^n\bi{n+k}k^2\eq& \f{4p^2E_{p-3}+(-1)^n}{4^{p-1}}
\\\eq& 4p^2E_{p-3}+(-1)^n(1-2p\,q_p(2)+3p^2q_p(2)^2)\ (\mo\ p^3).
\endaligned\tag4.16$$
For each $k=0,\ldots,n$, by (4.3) we have
$$\align &\bi{n+k}k^2
\\\eq&\f{\bi{2k}k^2}{16^k}\(1+p\sum_{j=1}^k\f1{2j-1}+\f{p^2}2\(\(\sum_{j=1}^k\f1{2j-1}\)^2
-\sum_{j=1}^k\f1{(2j-1)^2}\)\)^2
\\\eq&\f{\bi{2k}k^2}{16^k}\(1+2p\sum_{j=1}^k\f1{2j-1}+p^2\(2\(\sum_{j=1}^k\f1{2j-1}\)^2-\sum_{j=1}^k\f1{(2j-1)^2}\)\)
\\&\qquad\qquad\qquad\qquad\qquad\qquad\qquad\qquad\qquad\qquad\qquad\quad\pmod{p^3}.
\endalign$$
So we can obtain (4.16) by using (4.12)--(4.15).

Now we deduce (1.14). Combining (1.13) and (1.15) we get
$$\sum_{k=(p+1)/2}^{p-1}(21k+8)\bi{2k}k^3\eq(-1)^{(p+1)/2}32p^3E_{p-3}\pmod{p^4},$$
i.e.,
$$\sum_{k=1}^{(p-1)/2}(21(p-k)+8)\f{\bi{2(p-k)}{p-k}^3}{p^3}\eq(-1)^{(p+1)/2}32E_{p-3}\pmod p.$$
By [Su3, Lemma 2.1], for each $k=1,\ldots,(p-1)/2$ we have
$$\f{\bi{2(p-k)}{p-k}}p\eq\f{-2}{k\bi{2k}k}\pmod{p}.$$
Therefore
$$\sum_{k=1}^{p-1}(-21k+8)\l(\f{-2}{k\bi{2k}k}\r)^3\eq(-1)^{(p+1)/2}32E_{p-3}\pmod p,$$
which gives (1.14).

The proof of Theorem 1.4 is now complete. \qed

\heading{5. Some related conjectures}\endheading

We first pose the following conjecture similar to (1.6).

\proclaim{Conjecture 5.1} For any prime $p>3$ we have
$$\align\sum_{k=1}^{p-1}\f{\bi{2k}kH^{(2)}_k}k\eq&\f{2}3\cdot\f{H_{p-1}}{p^2}+\f{76}{135}p^2B_{p-5}\ (\mo\ p^3),
\\\sum_{k=1}^{p-1}\f{\bi{2k}kH_k^{(2)}}{k2^k}\eq&-\f3{16}\cdot\f{H_{p-1}}{p^2}+\f{479}{1280}p^2B_{p-5}\ (\mo\ p^3),
\\\sum_{k=1}^{p-1}\f{\bi{2k}kH_k^{(2)}}{k3^k}\eq&-\f 89\cdot\f{H_{p-1}}{p^2}+\f{268}{1215}p^2B_{p-5}\ (\mo\ p^3).
\endalign$$
\endproclaim
\Remark\ 5.1. It is known that
$$\f{H_{p-1}}{p^2}\eq -\f{B_{p-3}}3\pmod p\quad\t{for any prime}\ p>3$$
(see, e.g., [S1]).
\medskip

The following conjecture is close to Theorem 1.3.

\proclaim{Conjecture 5.2} Let $p>3$ be a prime. Then
$$\align\sum_{k=1}^{p-1}\f{\bi{2k}k}{k4^k}H_k\eq&\f 76pB_{p-3}\ (\mo\ p^2),
\\\sum_{k=1}^{p-1}\f{\bi{2k}k}{k4^k}H_{2k}\eq&\f 73pB_{p-3}\ (\mo\ p^2),
\\\sum_{k=1}^{p-1}\f{4^kH_{k-1}}{k^2\bi{2k}k}\eq&\f 23B_{p-3}\ (\mo\ p),
\\\sum_{k=1}^{(p-1)/2}\f{4^kH_{2k-1}}{k^2\bi{2k}k}\eq&\f 72 B_{p-3}\ (\mo\ p),
\\\sum_{k=1}^{p-1}\f{\bi{2k}k}{k4^k}H_k^{(2)}\eq&-\f 32\cdot\f{H_{p-1}}{p^2}+\f 7{80}p^2B_{p-5}\ (\mo\ p^3).
\endalign$$
Also,
$$\align\sum_{k=1}^{p-1}\f{\bi{2k}k}{k^24^k}H_k\eq&\f 32B_{p-3}\ (\mo\ p),
\\\sum_{k=1}^{p-1}\f{\bi{2k}k}{k^24^k}H_{2k}\eq&\f 52B_{p-3}\ (\mo\ p),
\endalign$$
and
$$\sum_{k=1}^{p-1}\f{\bi{2k}k}{k^24^k}\eq-\f{H_{(p-1)/2}^2}2-\f 74\cdot\f{H_{p-1}}{p}\ (\mo\ p^3)
\ \ \t{provided}\ p>5.$$
\endproclaim
\Remark\ 5.2. The author ever conjectured that
$$\sum_{k=1}^{(p-1)/2}\f{\bi{2k}k}{k4^k}H_{2k}\eq-2\l(\f{-1}p\r)E_{p-3}\ (\mo\ p)$$
for any prime $p>3$; this has been confirmed by his former student Hui-Qin Cao.
Using {\tt Mathematica 7} the author found that
$$\sum_{k=1}^\infty\f{4^kH_{k-1}}{k^2\bi{2k}k}=7\zeta(3),\ \ \sum_{k=1}^\infty\f{4^kH_{2k-1}}{k^2\bi{2k}k}=\f{21}2\zeta(3),$$
$$\sum_{k=1}^\infty\f{\bi{2k}k}{k4^k}H_k^{(2)}=\f32\zeta(3),\ \ \sum_{k=1}^\infty\f{\bi{2k}k}{k^24^k}=\f{\pi^2-3\log^24}6.$$
\medskip

Motivated by Theorem 4.1, we pose the following conjecture.
\proclaim{Conjecture 5.3} Let $p>3$ be a prime. Then
$$\align\sum_{k=1}^{p-1}\f{\bi{2k}k^2}{k16^k}H_{2k}^{(2)}\eq& B_{p-3}\ (\mo\ p),
\\\sum_{k=1}^{(p-1)/2}\f{\bi{2k}k^2}{k16^k}H_{2k}^{(2)}\eq&-\f 52B_{p-3}\ (\mo\ p),
\\\sum_{k=1}^{p-1}\f{\bi{2k}k^2}{k16^k}H_k^{(2)}\eq&-12\f{H_{p-1}}{p^2}+\f{7}{10}p^2B_{p-5}\ (\mo\ p^3),
\\\sum_{k=(p+1)/2}^{p-1}\f{\bi{2k}k^2}{k16^k}H_k^{(2)}\eq&\f{31}2p^2B_{p-5}\ (\mo\ p^3).
\endalign$$
Also,
$$\sum_{k=0}^{(p-3)/2}\f{\bi{2k}k^2}{(2k+1)16^k}H_k^{(2)}\eq-7B_{p-3}\pmod p$$
and
$$\sum_{k=(p+1)/2}^{p-1}\f{\bi{2k}k^2}{(2k+1)16^k}H_k^{(2)}\eq-\f{31}2p^2B_{p-5}\pmod {p^3}.$$
\endproclaim

\enddocument